\documentclass[a4paper,10pt]{article}

\usepackage{amsmath}
\usepackage{amssymb}
\usepackage{amsthm}
\usepackage{mathrsfs}
\usepackage{color}
\usepackage[dvips]{graphicx}
\usepackage[all]{xy}

\theoremstyle{plain}
\newtheorem{them}{Theorem}[section]
\newtheorem{lemma}[them]{Lemma}
\newtheorem{prop}[them]{Proposition}
\newtheorem{coro}[them]{Corollary}
\theoremstyle{definition}

\newtheorem{conj}[them]{Conjecture}

\newtheorem*{mthem}{Main Theorem}
\newtheorem*{conv}{Convention}

\begin{document}

\title{The zeta function of a finite category which has M\"obius inversion}
\author{Kazunori Noguchi \thanks{noguchi@math.shinshu-u.ac.jp}}
\date{}
\maketitle
\begin{abstract}
We prove certain conjecture holds true for a finite category which has M\"obius inversion. The conjecture states a relationship between the zeta function of a finite category and the Euler characteristic of a finite category.
\end{abstract}


\thispagestyle{empty}

\section{Introduction}

In \cite{NogB}, the zeta function of a finite category was defined and one conjecture was proposed. The zeta function of a finite category $I$ is defined by $$\zeta_I(z)=\exp\left( \sum_{m=1}^{\infty} \frac{\# N_m(I)}{m} z^m\right)$$
where $$N_m(I)=\{\xymatrix{(x_0\ar[r]^{f_1}&x_1\ar[r]^{f_2}&\dots\ar[r]^{f_m}&x_m)} \text{ in } I \}.$$  And the conjecture states a relationship between the zeta function of a finite category and the Euler characteristic of a finite category. 
\begin{conj}
Suppose $I$ is a finite category and its series Euler characteristic $\chi_{\Sigma}(I)$ exists. Then, we have
\renewcommand{\theenumi}{C\arabic{enumi}}
\renewcommand{\labelenumi}{(\theenumi)}
\begin{enumerate}
\item the zeta function of $I$ is a finite product of the following form $$\zeta_I(z)=\prod \frac{1}{(1-\alpha_iz)^{\beta_i}}\exp\left( \sum \frac{\gamma_j z^j}{j(1-\delta_j z)^j}\right)$$
for some complex numbers $\alpha_i, \beta_i, \gamma_j,\delta_j $.
\item  $\displaystyle \sum \beta_i$ is the number of objects of $I$ 
\item each $\alpha_i$ is an eigen value of $A_I$. Hence, $\alpha_i$ is an algebraic integer.
\item $\displaystyle \sum \frac{\beta_i}{\alpha_i} +\sum (-1)^j\frac{\gamma_j}{\delta_j^{j+1}}=\chi_{\Sigma}(I).$
\end{enumerate}
\end{conj}

In \cite{NogA}, it was verified the conjecture holds true in concrete cases, that is, when a finite category is a groupoid, an acyclic category, a category which has two objects and so on.

In this paper, we prove the conjecture holds true when a finite category has M\"obius inversion. A finite category $I$ has \textit{M\"obius inversion} if its adjacency matrix  $A_I$ has the inverse matrix where $A_I$ is an $N\times N$-matrix whose $(i,j)$-entry is the number of morphisms of $I$ from $x_i$ to $x_j$ when the set of objects of $I$ is $$\mathrm{Ob}(I)=\{x_1,x_2,\dots,x_N\}.$$ This class of finite categories is large and very important to consider the Euler characteristic of a finite category.

Now, let us recall the Euler characteristic of a category. It is defined in various ways, that is, we have Leinster's Euler characteristic $\chi_{L}$ \cite{Leia}, the series Euler characteristic $\chi_{\Sigma}$ by Berger-Leinster \cite{Leib}, the $L^2$-Euler characteristic $\chi^{(2)}$ by Fiore-L\"uck-Sauer \cite{FLS}, the extended $L^2$-Euler characteristic $\chi^{(2)}_{\mathrm{ex}}$ and the Euler characteristic for $\mathbb{N}$-filtered acyclic categories $\chi_{\mathrm{fil}}$ by the author \cite{NogA},\cite{NogB}. Here, we focus on the series Euler characteristic which is used in the conjecture. For a finite category $I$ whose set of objects is $\{x_1,\dots, x_N\}$, its series Euler characteristic $\chi_{\Sigma}(I)$ is defined by substituting $-1$ to $t$ of 
$$\frac{\mathrm{sum}(\mathrm{adj}(E-(A_I-E)t))}{\det(E-(A_I-E)t)}$$
if it exists where sum means to take the sum of all of the entries of a matrix. When $I$ has M\"obius inversion, $\chi_{\Sigma}(I)=\mathrm{sum}(A_I^{-1})$. So, in order to read this paper it is enough that we recognize the Euler characteristic of a finite category $I$ is defined by $\mathrm{sum}(A_I^{-1})$ since all of the categories in this paper have M\"obius inversion.

There are two reasons why we importantly deal with the class of finite categories which has M\"obius inversion. The one is that the two Euler characteristics $\chi_L$ and $\chi_{\Sigma}$ coincide for this class. The another is that important finite categories are contained in the class, for example, posets, acyclic categories and groups (with one-object). In particular, when a category is a finite acyclic category all of those Euler characteristics take the same value.

The following is our main theorem.

\begin{mthem}
Let $I$ be a finite category which has $N$-objects. Suppose $\det A_I \not= 0$ and the polynomial $\det(E-A_Iz)$ is factored to the following form:
$$\det(E-A_Iz)=(-1)^N \det A_I (z-\alpha_1)^{e_1}\dots(z-\alpha_n)^{e_n}$$ where each $e_i\ge 1$ and $\alpha_i\not = \alpha_j$ if $i\not = j$. And suppose the rational function $\frac{\mathrm{sum}(\mathrm{adj}(E-A_I z)A_I)}{\det (E-A_I z)}$ has a partial fraction decomposition to the following form: 
$$\frac{\mathrm{sum}(\mathrm{adj}(E-A_I z)A_I)}{\det (E-A_I z)}=\frac{1}{(-1)^N\det A_I}\sum^n_{k=1}\sum^{e_k}_{j=1}\frac{A_{k,j}}{(z-\alpha_k)^j}.$$ 
\begin{enumerate}
\item  Then the zeta function of $I$ is 
\begin{multline*}
\zeta_I(z)=\prod^n_{k=1}\frac{1}{(1-\frac{1}{\alpha_k}z)^{-\frac{A_{k,1}}{(-1)^N\det A_I}}}  \\ \times \exp\left( \frac{1}{(-1)^N\det A_I}  \sum^n_{k=1} \sum^{e_k -1}_{j=1} \frac{z^j}{j(1-\frac{1}{\alpha_k}z)^j}\bigg( \sum_{i=j}^{e_k-1} j \binom{i}{j} (-1)^{i+1} \bigg(\frac{1}{\alpha_k}\bigg)^{i+j} \frac{A_{k,i+1}}{i}\bigg)\right)
\end{multline*}
\item $\displaystyle \sum^n_{k=1}-\frac{A_{k,1}}{(-1)^N\det A_I}=N$
\item  Each $\frac{1}{\alpha_k}$ is an eigen value of $A_I$. In particular, $\frac{1}{\alpha_k}$ is an algebraic integer
\item  
\begin{multline*}
\displaystyle \sum^n_{k=1}\frac{-\frac{A_{k,1}}{(-1)^N\det A_I}}{\frac{1}{\alpha_k}}\\
+ \frac{1}{(-1)^N\det A_I} \sum^n_{k=1} \sum ^{e_k-1}_{j=1} (-1)^j\frac{ \sum_{i=j}^{e_k-1}  \binom{i-1}{j-1} (-1)^{i-1} \bigg(\frac{1}{\alpha_k}\bigg)^{i+j} A_{k,i+1}}{\left(\frac{1}{\alpha_k}\right)^{j+1}}=\chi_{\Sigma}(I).
\end{multline*}
\end{enumerate}
\end{mthem}

Consequently, this theorem implies a corollary which states about behavior of singular points and zeros of the zeta function of a finite category. Namely, when $z$ is a complex number, $z$ is a singular point or zero of $\zeta_I$ if and only if $z$ is a root of $\det(E-A_Iz)$. Studying about singular points and zeros of a zeta function is very important problem for almost zeta functions.

This paper is organized as follows. In section \ref{prep}, we prepare a proposition and some lemmas to prove our main theorem. In section \ref{main}, we give a proof of our main theorem.

\section{Preparation}\label{prep}
\begin{conv}
Throughout this paper, $I$ is a finite category which has $N$-objects $$\mathrm{Ob}(I)=\{x_1,x_2,\dots,x_N\}.$$ And $A_I$ is an $N\times N$-matrix whose $(i,j)$-entry is the number of morphisms from $x_i$ to $x_j$ in $I$.
\end{conv}
The proof of our main theorem will start from the following integral expression:

\begin{prop}\label{integral}
 If $\det A_I\not = 0$, the zeta function of $I$ is
 $$\zeta_I(z)=\exp\left(\int \frac{\mathrm{sum}(\mathrm{adj}(E-A_I z)A_I)}{\det (E-A_I z)}dz +C \right)$$
 for some constant $C$.
 \begin{proof}
 We have 
 \begin{eqnarray*}
 z\frac{d}{dz} \log \zeta_I(z)&=& z\frac{d}{dz}\left(\sum^{\infty}_{m=1} \frac{\# N_m}{m}z^m \right) \\
 &=&\sum^{\infty}_{m=1}\# N_m z^m \\
 &=&\sum^{\infty}_{m=1}\mathrm{sum}(A_I^m) z^m \\
 &=&\mathrm{sum}\bigg( \sum^{\infty}_{m=1}(A_I^m) z^m \bigg). 
 \end{eqnarray*}
Put $F(z)=\sum^{\infty}_{m=1}(A_I^m) z^m$. Then, we have
\begin{eqnarray*}
\mathrm{adj}(E-A_I z) (E-A_I z)&=&\det (E-A_I z)E \\
\mathrm{adj}(E-A_I z) (E-A_I z)F(z)&=&\det (E-A_I z)F(z) \\
\mathrm{adj}(E-A_I z)A_I z&=&\det (E-A_I z)F(z). 
\end{eqnarray*}
The polynomial $\det(E-A_I z)$ is not the zero polynomial since $\det(E-A_I z)=1$ for $z=0$. Hence, we have
$$F(z)=\frac{\mathrm{adj}(E-A_I z)A_I z}{\det (E-A_I z)}.$$
Therefore, we have $$\mathrm{sum}(F(z))=\frac{\mathrm{sum}(\mathrm{adj}(E-A_I z)A_I z)}{\det (E-A_I z)}.$$
Thus, we obtain
\begin{eqnarray*}
\zeta_I(z)&=&\exp \left( \int \frac{1}{z} \mathrm{sum}(F(z))dz +C \right) \\
&=&\exp\left(\int \frac{\mathrm{sum}(\mathrm{adj}(E-A_I z)A_I)}{\det (E-A_I z)}dz +C \right).
\end{eqnarray*}
 \end{proof}
\end{prop}

\begin{lemma}\label{lemma1}
We have $$\mathrm{sum}(\mathrm{adj}(E-A_Iz)A_I)=\frac{1}{z}\bigg(\mathrm{sum}(\mathrm{adj}(E-A_Iz))-N\det (E-A_Iz)\bigg).$$
\begin{proof}
\begin{eqnarray*}
\mathrm{adj}(E-A_Iz)(E-A_Iz)&=&\det (E-A_Iz)E \\
\mathrm{adj}(E-A_Iz)-\mathrm{adj}(E-A_Iz)A_Iz&=&\det (E-A_Iz)E\\
\mathrm{adj}(E-A_Iz)A_I&=&\frac{1}{z}\bigg(\mathrm{adj}(E-A_Iz)-\det (E-A_Iz)E\bigg) \\
\mathrm{sum}(\mathrm{adj}(E-A_Iz)A_I)&=&\frac{1}{z}\bigg(\mathrm{sum}(\mathrm{adj}(E-A_Iz)) \\
&&-N\det (E-A_Iz)\bigg).
\end{eqnarray*}
\end{proof}
\end{lemma}

\begin{lemma}\label{lemmad}
Let $$\det (E-A_Iz)=d_0+d_1z+\dots+d_{N-1}z^{N-1}+d_Nz^N.$$
Then, $d_0=1$ and $d_N=(-1)^N \det A_I$.
\begin{proof}
For $z=0$, we have $$1=\det (E-A_I 0)=d_0.$$ Now, we express the unit matrix $E$ by Kronecker's delta, that is, $E=(\delta_{ij})_{ij}$ and let $A_I=(a_{ij})_{ij}$. Then, we have
\begin{eqnarray*}
\det(E-A_Iz)&=&\det\big( (\delta_{ij}-a_{ij}z)_{ij}\big)\\
&=&\sum_{\sigma \in S_N} \mathrm{sign} (\sigma) (\delta_{1\sigma(1)}-a_{1\sigma(1)}z)\dots (\delta_{N\sigma(N)}-a_{N\sigma(N)}z).
\end{eqnarray*}
Hence, the coefficient of $z^N$ in $\det(E-A_Iz)$ is $$\sum_{\sigma \in S_N} \mathrm{sign} (\sigma) (-a_{1\sigma(1)})\dots (-a_{N\sigma(N)}).$$ Therefore, we obtain $d_N=(-1)^N\det A_I$.
\end{proof}
\end{lemma}

\begin{lemma}\label{lemmak}
Let $$\mathrm{sum}\big(\mathrm{adj}(E-A_Iz)\big)=k_0+k_1z+\dots +k_{N-1}z^{N-1}.$$
Then, $k_0=N$ and $k_{N-1}=(-1)^{N-1}\mathrm{sum}(\mathrm{adj}(A_I))$.
\begin{proof}
Let $A_{ij},A'_{ij}$ be the $(i,j)$-cofactor of $A_I$ and $E-A_Iz$, respectively. If $i=j$, the same argument of Lemma \ref{lemmad} implies that $A'_{ii}$ is a polynomial whose constant term is 1 and term of degree $N-1$ is $(-1)^{N-1}A_{ii}z^{N-1}$. If $i\not = j$, $A'_{ij}$ is a polynomial whose constant term is 0 and term of degree $N-1$ is $(-1)^{N-1}A_{ij}z^{N-1}$. Hence, we obtain the result.
\end{proof}
\end{lemma}

\begin{lemma}\label{lemmam}
Let $$\mathrm{sum}(\mathrm{adj}(E-A_I z)A_I)=m_0+m_1z+\dots+m_{N-2}z^{N-2}+m_{N-1}z^{N-1}.$$
Then, $m_{N-1}=(-1)^{N+1}N\det A_I$ and $m_{N-2}=(-1)^{N-1}\mathrm{sum}(\mathrm{adj}(A_I))-N d_{N-1}$
where $d_{N-1}$ is the same coefficient defined in Lemma \ref{lemmad}.
\begin{proof}
Lemma \ref{lemma1}, Lemma \ref{lemmad} and Lemma \ref{lemmak} imply this result.
\end{proof}
\end{lemma}

\begin{lemma}\label{lemmadeg}
If $\det A_I\not =0$, then the degree of $\mathrm{sum}(\mathrm{adj}(E-A_Iz)A_I)$ is less than the degree of $\det (E-A_Iz).$
\begin{proof}
Lemma \ref{lemmad} implies the degree of the polynomial $\det(E-A_Iz)$ is $N$ since $\det A_I\not=0$. And the degree of the polynomial $\mathrm{sum}(\mathrm{adj}(E-A_Iz))$ is less than or equal to $N-1$. Hence, Lemma \ref{lemma1} implies 
\begin{eqnarray*}
\deg\big( \mathrm{sum}(\mathrm{adj}(E-A_Iz)A_I) \big)&=&\deg \bigg( \frac{1}{z}\bigg(\mathrm{sum}(\mathrm{adj}(E-A_Iz)) \\ 
&& -N\det (E-A_Iz)\bigg)\bigg)\\
&=&\deg \big( -N\det (E-A_Iz) \big) -1 \\
&=&N-1.
\end{eqnarray*}
\end{proof}
\end{lemma}

\begin{lemma}\label{lemmack}
Suppose $\alpha$ is a non-zero complex number and $n$ is a natural number. Then, we have
$$\frac{1}{(z-\alpha)^n}=\frac{C_1z}{(z-\alpha)}+\frac{C_2z^2}{(z-\alpha)^2}+\dots+\frac{C_n z^n}{(z-\alpha)^n}+\left(-\frac{1}{\alpha}\right)^n$$
where each $C_k$ is $\binom{n}{k}(-\frac{1}{\alpha})^n(-1)^k$.
\begin{proof}
We have
\begin{eqnarray*}
\frac{1}{(z-\alpha)^n}&=&\frac{1}{(z-\alpha)^n}-\left(-\frac{1}{\alpha}\right)^n+\left(-\frac{1}{\alpha}\right)^n\\
&=&\frac{-\sum^n_{k=1}\binom{n}{k}(-\frac{z}{\alpha})^k }{(z-\alpha)^n}+\left(-\frac{1}{\alpha}\right)^n
\end{eqnarray*}
Here, put $$\frac{-\sum^n_{k=1}\binom{n}{k}(-\frac{z}{\alpha})^k }{(z-\alpha)^n}=\frac{C_1z}{(z-\alpha)}+\frac{C_2z^2}{(z-\alpha)^2}+\dots+\frac{C_n z^n}{(z-\alpha)^n}.$$
Then, we show each $C_k$ is $\binom{n}{k}(-\frac{1}{\alpha})^n(-1)^k$ by induction on $k$. By a reduction to common denominator of the right hand side, we have 

\begin{eqnarray}\frac{-\sum^n_{k=1}\binom{n}{k}(-\frac{z}{\alpha})^k }{(z-\alpha)^n}&=&\frac{C_1z(z-\alpha)^{n-1}+C_2z^2(z-\alpha)^{n-2}+\dots+C_n z^n}{(z-\alpha)^n}.\label{eq1}
\end{eqnarray}
By observing the numerators of both side, we have $C_1z(-\alpha)^{n-1}=-\binom{n}{1}(-\frac{z}{\alpha})$. Hence, $C_1=-\frac{n}{(-\alpha)^n}$. Next, we assume that it holds for any $i\le k$. By comparing the terms of degree $k+1$ of both sides of (\ref{eq1}), we have
\begin{eqnarray*}
-\binom{n}{k+1}\bigg(-\frac{z}{\alpha}\bigg)^{k+1}&=&\sum_{j=1}^{k+1} \binom{n-j}{k+1-j}z^{k+1}(-\alpha)^{n-(k+1)}C_j \\
&=&\sum_{j=1}^{k} \binom{n-j}{k+1-j}z^{k+1}(-\alpha)^{n-(k+1)}C_j + \\
&& z^{k+1}(-\alpha)^{n-(k+1)}C_{k+1} \\
&=&\sum_{j=1}^{k} \binom{n-j}{k+1-j}z^{k+1}(-\alpha)^{n-(k+1)}\binom{n}{j}\bigg(-\frac{1}{\alpha}\bigg)^n(-1)^j \\ 
&&+z^{k+1}(-\alpha)^{n-(k+1)}C_{k+1} \\
&=&\sum_{j=1}^{k} \binom{k+1}{j}\binom{n}{k+1}z^{k+1}\bigg(-\frac{1}{\alpha}\bigg)^{k+1}(-1)^j \\
&&+z^{k+1}(-\alpha)^{n-(k+1)}C_{k+1}.
\end{eqnarray*}
Hence, we have
\begin{eqnarray*}
-\binom{n}{k+1}&=&\sum^k_{j=1} \binom{k+1}{j}\binom{n}{k+1}(-1)^j +C_{k+1}(-\alpha)^n \\
C_{k+1}(-\alpha)^n &=&-\binom{n}{k+1}\bigg( 1+\sum^k_{j=1} \binom{k+1}{j}(-1)^j \bigg) \\
C_{k+1}(-\alpha)^n &=&-\binom{n}{k+1}\bigg( \sum^k_{j=0} \binom{k+1}{j}(-1)^j \bigg) \\
C_{k+1}(-\alpha)^n &=&-\binom{n}{k+1}(-1)^{k+1} \\
C_{k+1}&=&\binom{n}{k+1}\bigg( -\frac{1}{\alpha}\bigg)^n(-1)^{k+1}.
\end{eqnarray*}
Hence, we obtain the result.
\end{proof}
\end{lemma}


\section{The proof of main theorem}\label{main}

\begin{them}\label{noguchi}
Suppose $\det A_I \not= 0$ and the polynomial $\det(E-A_Iz)$ is factored to the following form:
$$\det(E-A_Iz)=(-1)^N \det A_I (z-\alpha_1)^{e_1}\dots(z-\alpha_n)^{e_n}$$ where each $e_i\ge 1$ and $\alpha_i\not = \alpha_j$ if $i\not = j$. And suppose the rational function $\frac{\mathrm{sum}(\mathrm{adj}(E-A_I z)A_I)}{\det (E-A_I z)}$ has a partial fraction decomposition to the following form: 
$$\frac{\mathrm{sum}(\mathrm{adj}(E-A_I z)A_I)}{\det (E-A_I z)}=\frac{1}{(-1)^N\det A_I}\sum^n_{k=1}\sum^{e_k}_{j=1}\frac{A_{k,j}}{(z-\alpha_k)^j}.$$ 
\begin{enumerate}
\item \label{C1} Then the zeta function of $I$ is 
\begin{multline*}
\zeta_I(z)=\prod^n_{k=1}\frac{1}{(1-\frac{1}{\alpha_k}z)^{-\frac{A_{k,1}}{(-1)^N\det A_I}}}  \\ \times \exp\left( \frac{1}{(-1)^N\det A_I}  \sum^n_{k=1} \sum^{e_k -1}_{j=1} \frac{z^j}{j(1-\frac{1}{\alpha_k}z)^j}\bigg( \sum_{i=j}^{e_k-1}  \binom{i-1}{j-1} (-1)^{i-1} \bigg(\frac{1}{\alpha_k}\bigg)^{i+j} A_{k,i+1}\bigg)\right)
\end{multline*}
\item \label{C2}$\displaystyle \sum^n_{k=1}-\frac{A_{k,1}}{(-1)^N\det A_I}=N$
\item \label{C3} Each $\frac{1}{\alpha_k}$ is an eigen value of $A_I$. In particular, $\frac{1}{\alpha_k}$ is an algebraic integer
\item \label{C4} 
\begin{multline*}
\displaystyle \sum^n_{k=1}\frac{-\frac{A_{k,1}}{(-1)^N\det A_I}}{\frac{1}{\alpha_k}}\\
+ \frac{1}{(-1)^N\det A_I} \sum^n_{k=1} \sum ^{e_k-1}_{j=1} (-1)^j\frac{ \sum_{i=j}^{e_k-1}  \binom{i-1}{j-1} (-1)^{i-1} \bigg(\frac{1}{\alpha_k}\bigg)^{i+j} A_{k,i+1}}{\left(\frac{1}{\alpha_k}\right)^{j+1}}=\chi_{\Sigma}(I).
\end{multline*}
\end{enumerate}

\begin{proof}
We first show (\ref{C1}). Lemma \ref{lemmadeg} implies the degree of the numerator is less than the degree of the denominator for the rational function $$\frac{\mathrm{sum}(\mathrm{adj}(E-A_I z)A_I)}{\det (E-A_I z)}.$$
Hence, we can have a partial fraction decomposition of the following form $$\frac{\mathrm{sum}(\mathrm{adj}(E-A_I z)A_I)}{\det (E-A_I z)}=\frac{1}{(-1)^N\det A_I}\sum^n_{k=1}\sum^{e_k}_{j=1}\frac{A_{k,j}}{(z-\alpha_k)^j}$$for some complex numbers $A_{k,j}$.
Proposition \ref{integral} implies
\begin{eqnarray*}
\zeta_I(z)&=&\exp\left( \int \frac{1}{(-1)^N\det A_I}\sum^n_{k=1}\sum^{e_k}_{j=1}\frac{A_{k,j}}{(z-\alpha_k)^j} dz +C\right)\\
&=&\exp\bigg( \frac{1}{(-1)^N\det A_I}\int \sum^n_{k=1}\frac{A_{k,1}}{(z-\alpha_k)} dz +\\
&& \frac{1}{(-1)^N\det A_I} \int \sum^n_{k=1}\sum^{e_k}_{j=2}\frac{A_{k,j}}{(z-\alpha_k)^j} dz+C\bigg) \\
&=&\exp\bigg( \frac{1}{(-1)^N\det A_I} \sum^n_{k=1} A_{k,1}\log (z-\alpha_k)  +\\
&& \frac{1}{(-1)^N\det A_I} \sum^n_{k=1}\sum^{e_k}_{j=2}-\frac{A_{k,j}}{(j-1)}\frac{1}{(z-\alpha_k)^{j-1}} +C\bigg)\\
&=&\prod_{k=1}^n \frac{1}{(z-\alpha_k)^{-\frac{A_{k,1}}{(-1)^N\det A_I}}}\times\\
&&\exp\bigg(\frac{1}{(-1)^N\det A_I} \sum^n_{k=1}\sum^{e_k-1}_{j=1}-\frac{A_{k,j+1}}{j}\frac{1}{(z-\alpha_k)^{j}} \bigg)\exp{C}\\
&=&\prod_{k=1}^n \frac{1}{(-\alpha_k)^{-\frac{A_{k,1}}{(-1)^N\det A_I}}(1-\frac{1}{\alpha_k}z)^{-\frac{A_{k,1}}{(-1)^N\det A_I}}}\times\\
&&\exp\bigg(\frac{1}{(-1)^N\det A_I} \sum^n_{k=1}\sum^{e_k-1}_{j=1}-\frac{A_{k,j+1}}{j}\frac{1}{(z-\alpha_k)^{j}} \bigg)C'\\
&=&\prod_{k=1}^n \frac{1}{(1-\frac{1}{\alpha_k}z)^{-\frac{A_{k,1}}{(-1)^N\det A_I}}}\times\\
&&\exp\bigg(\frac{-1}{(-1)^N\det A_I} \sum^n_{k=1}\sum^{e_k-1}_{j=1}\frac{A_{k,j+1}}{j}\frac{1}{(z-\alpha_k)^{j}} \bigg)C''
\end{eqnarray*}
where we did and will replace the constant term as $C,C',C''\dots$. Lemma \ref{lemmack} implies
\begin{eqnarray*}
\zeta_I(z)&=&\prod_{k=1}^n \frac{1}{(1-\frac{1}{\alpha_k}z)^{-\frac{A_{k,1}}{(-1)^N\det A_I}}}\times\\
&&\exp\bigg(\frac{-1}{(-1)^N\det A_I} \sum^n_{k=1}\sum^{e_k-1}_{j=1}\frac{A_{k,j+1}}{j} \sum_{i=1}^j \frac{\binom{j}{i}(-\frac{1}{\alpha_k})^j (-z)^i}{(z-\alpha_k)^i}\bigg)C'''
\end{eqnarray*}
Here, we use the boundary condition $\zeta_{I}(0)=1$. This condition is directly implied by the definition of the zeta function. Hence, we obtain $C'''=1$. By exchanging $\sum_i$ and $\sum_j$, we have
\begin{multline*}
\zeta_I(z)=\prod^n_{k=1}\frac{1}{(1-\frac{1}{\alpha_k}z)^{-\frac{A_{k,1}}{(-1)^N\det A_I}}}  \\ \times \exp\left( \frac{1}{(-1)^N\det A_I}  \sum^n_{k=1} \sum^{e_k -1}_{j=1} \frac{z^j}{j(1-\frac{1}{\alpha_k}z)^j}\bigg( \sum_{i=j}^{e_k-1} j \binom{i}{j} (-1)^{i+1} \bigg(\frac{1}{\alpha_k}\bigg)^{i+j} \frac{A_{k,i+1}}{i}\bigg)\right)
\end{multline*}
Hence, we obtain the result.

Next we show (\ref{C2}). We observe the numerators of both sides $$\frac{\mathrm{sum}(\mathrm{adj}(E-A_I z)A_I)}{\det (E-A_I z)}=\frac{1}{(-1)^N\det A_I} \sum^n_{k=1}\sum^{e_k}_{j=1}\frac{A_{k,j}}{(z-\alpha_k)^j}.$$
The numerator of the left hand side is $m_0+\dots +m_{N-1}z^{N-1}$ and Lemma \ref{lemmam} implies $m_{N-1}=(-1)^{N+1}N\det A_I$. For the right hand side, when it is transformed to the left hand side by a reduction to common denominator, the coefficient of $z^{N-1}$ of the numerator is $\sum_{k=1}^n A_{k,1}$. Hence, we have $\sum_{k=1}^n A_{k,1}=(-1)^{N+1}N\det A_I$. Thus, we obtain 
$$\sum^n_{k=1}-\frac{A_{k,1}}{(-1)^N{det A_I}}=N.$$

 We show (\ref{C3}). Since each $\alpha_k$ is a root of the polynomial $\det(E-A_Iz)$, we obtain
 \begin{eqnarray*}
 \det(E-A_I\alpha_k)&=&0 \\
 (\alpha_k)^N \det(E\frac{1}{\alpha_k}-A_I)&=&0.
 \end{eqnarray*}
 Hence, $\frac{1}{\alpha_k}$ is an eigen value of $A_I$. Note that $\alpha_k\not =0$. Moreover, since $\det(E\lambda-A_I)$ is a monic polynomial with coefficients in $\mathbb{Z}$, $\frac{1}{\alpha_k}$ is an algebraic integer.

 Finally, we show (\ref{C4}). By $(\ref{C1})$, we have
 \begin{multline*}
\zeta_I(z)=\prod^n_{k=1}\frac{1}{(1-\frac{1}{\alpha_k}z)^{-\frac{A_{k,1}}{(-1)^N\det A_I}}}  \\ \times \exp\left( \frac{1}{(-1)^N\det A_I}  \sum^n_{k=1} \sum^{e_k -1}_{j=1} \frac{z^j}{j(1-\frac{1}{\alpha_k}z)^j}\bigg( \sum_{i=j}^{e_k-1} j \binom{i}{j} (-1)^{i+1} \bigg(\frac{1}{\alpha_k}\bigg)^{i+j} \frac{A_{k,i+1}}{i}\bigg)\right)
\end{multline*}
The left hand side of (C\ref{C4}) is
\begin{eqnarray*}
(C\ref{C4})&=&\displaystyle \sum^n_{k=1}\frac{-\frac{A_{k,1}}{(-1)^N\det A_I}}{\frac{1}{\alpha_k}} \\
&&+ \frac{1}{(-1)^N\det A_I} \sum^n_{k=1} \sum ^{e_k-1}_{j=1} (-1)^j\frac{ \sum_{i=j}^{e_k-1}  \binom{i-1}{j-1} (-1)^{i-1} \bigg(\frac{1}{\alpha_k}\bigg)^{i+j} A_{k,i+1}}{\left(\frac{1}{\alpha_k}\right)^{j+1}} \\
&=&\sum^n_{k=1}\bigg(-\frac{\alpha_k A_{k,1}}{(-1)^N\det A_I} \\
&&+\frac{1}{(-1)^N \det A_I} \sum^{e_k-1}_{j=1}\sum^{e_k-1}_{i=j} (-1)^{j}\bigg(-\frac{1}{\alpha_k}\bigg)^{i-1}\binom{i-1}{j-1}A_{k,i+1} \\
&=&\sum^n_{k=1}\bigg(-\frac{\alpha_k A_{k,1}}{(-1)^N\det A_I} \\
&&+\frac{1}{(-1)^N\det A_I} \sum^{e_k-1}_{i=1} \bigg(-\frac{1}{\alpha_k}\bigg)^{i-1}A_{k,i+1}\big(\sum^i_{j=1}(-1)^j\binom{i-1}{j-1}\big)\\
&=&\frac{1}{(-1)^N\det A_I} \bigg( \sum^n_{k=1}(-\alpha_k A_{k,1}-A_{k,2})\bigg).
\end{eqnarray*}
 So it is enough to show $$\frac{1}{(-1)^N\det A_I}\bigg( \sum^n_{k=1}-\alpha_kA_{k,1}-A_{k,2}\bigg)=\mathrm{sum}(A_I^{-1}).$$  Therefore, it suffices to show $$\sum^n_{k=1}(-\alpha_kA_{k,1}-A_{k,2})=(-1)^N \mathrm{sum(\mathrm{adj}(A_I))}.$$
 We have $$\frac{1}{(-1)^N\det A_I}\bigg(\frac{A_{1,1}}{(z-\alpha_1)}+\dots+\frac{A_{n,e_n}}{(z-\alpha_n)^{e_n}}\bigg)= \frac{\mathrm{sum}(\mathrm{adj}(E-A_I z)A_I)}{\det (E-A_I z)}.$$
 We have already shown $\sum_{k=1}^n A_{k,1}=(-1)^{N+1}N\det A_I$.
 Now, we observe the term of degree $N-2$ of both sides. Then, Lemma \ref{lemmam} implies
\begin{multline} 
\sum^n_{k=1}A_{k,2}+\sum^n_{k=1}A_{k,1}(-\alpha_1 e_1-\dots-\alpha_k(e_k-1)\dots -\alpha_ne_n)=m_{N-2}= \\(-1)^{N-1}\mathrm{sum}(\mathrm{adj}A_I)-Nd_{N-1}.
\end{multline}
Here, we observe the equation
\begin{eqnarray*}
\det(E-A_Iz)&=&(-1)^N\det A_I (z-\alpha_1)^{e_1}\dots (z-\alpha_n)^{e_n} \\
&=&d_0+d_1z+\dots+d_{N-1}z^{N-1}+z^N.
\end{eqnarray*}
We obtain $d_{N-1}=(-1)^N \det A_I (\sum_{k=1}^n -\alpha_k e_k)$.
Therefore, 
\begin{multline}
\sum_{k=1}^n A_{k,2}=(-1)^{N-1}\mathrm{sum}(\mathrm{adj}A_I)-N(-1)^N\det A_I (\sum_{k=1}^n-\alpha_k e_k) -\\
\sum_{k=1}^nA_{k,1}((-\alpha_1e_1-\dots-\alpha_k(e_k-1)\dots -\alpha_ne_n)).
\end{multline}
Hence, we obtain
\begin{eqnarray*}
\sum^n_{k=1}(-\alpha_k A_{k,1}-A_{k,2})&=&\sum_{k=1}^n A_{k,1}(\sum_{k=1}^n -\alpha_k e_k)+\\
&&(-1)^{N}\mathrm{sum}(\mathrm{adj}A_I)+N(-1)^N\det A_I (\sum_{k=1}^n -\alpha_k e_k) \\
&=&\bigg(\sum_{k=1}^n -\alpha_k e_k\bigg)\bigg(\sum_{k=1}^n A_{k,1}+N(-1)^N\det A_I\bigg) +\\
&&(-1)^{N}\mathrm{sum}(\mathrm{adj}(A_I)) \\
&=&(-1)^{N}\mathrm{sum}(\mathrm{adj}(A_I)).
\end{eqnarray*}
Thus, we obtain the result.
\end{proof}
\end{them}

It is very important to study about behavior of singular points and zeros of a zeta function. The following corollary gives us a little information about the problem.

\begin{coro}
Suppose $\det A_I\not =0$ and $z$ is a complex number. Then, $z$ is a singular point or zero of $\zeta_I$ if and only if $z$ is a root of $\det (E-A_Iz)$.
\begin{proof}
Theorem \ref{noguchi} (\ref{C1}) directly implies all of the singular points and zeros are roots of $\det (E-A_Iz)$. Conversely, suppose $z$ is a root of $\det (E-A_Iz)$ but $z$ is not a singular point and zero. Then, $z=\alpha_l$ for some $k$. The index $-\frac{A_{l,1}}{(-1)^N \det A_I}$ must be 0. Namely, we have $A_{l,1}=0$. For $j=e_l-1$, $\sum^{e_l-1}_{i=e_l-1}-\frac{A_{l,i+1}}{i}\binom{i}{j}(-\frac{1}{\alpha_l})^{i-(e_l-1)}$ must be 0 since $\zeta_I(z)$ is defined. Hence, we have $A_{l,e_l}=0$. As this, we can show each $A_{l,j}=0$ by the descent from $j=e_l-1$. Hence, we have 
\begin{eqnarray*}
\frac{\mathrm{sum}(\mathrm{adj}(E-A_I z)A_I)}{\det (E-A_I z)}&=&\frac{1}{(-1)^N\det A_I}\sum^n_{k=1}\sum^{e_l}_{j=1}\frac{A_{k,j}}{(z-\alpha_l)^j}\\
&=&\frac{1}{(-1)^N\det A_I}\sum^n_{k=1,k\not=l}\sum^{e_l}_{j=1}\frac{A_{k,j}}{(z-\alpha_l)^j}
\end{eqnarray*}
Hence, we obtain
\begin{eqnarray*}
\det(E-A_I z)&=&(-1)^N\det A_I (z-\alpha_1)^{e_1}\dots(z-\alpha_n)^{e_n} \\
&=&(-1)^N\det A_I (z-\alpha_1)^{e_1}\dots \\
&&(z-\alpha_{l-1})^{e_{l-1}}(z-\alpha_{l+1})^{e_{l+1}}\dots(z-\alpha_n)^{e_n} 
\end{eqnarray*}
The polynomial $\det(E-A_I z)$ has two different degrees since each $e_k\ge 1$. This contradiction implies $z=\alpha_l$ is a singular point or zero of $\zeta_I$. 
\end{proof}
\end{coro}


\begin{thebibliography}{AAA99}

\bibitem [BL08]{Leib}  C. Berger and T. Leinster. The Euler characteristic of a category as the sum of a divergent series, \textit{Homology, Homotopy Appl., }10(1):41-51, 2008.



\bibitem[FLS11]{FLS} T. M. Fiore, W. L\"{u}ck and R. Sauer. Finiteness obstructions and Euler characteristics of categories, \textit{Adv. Math}, Vol. 226, Number 3, (2011), 2371--2469.



\bibitem [Lei08]{Leia}    T. Leinster. The Euler characteristic of a category,\textit{ Doc. Math.}, 13:21-49, 2008, arXiv:math.CT/0610260




\bibitem [Nog11]{Nog11} K. Noguchi. The Euler characteristic of acyclic categories. Kyushu Journal of Mathematics, vol. 65 No.1 (2011), 85-99.

\bibitem [NogA]{NogA} K. Noguchi. The Euler characteristics of categories and the barycentric subdivision. arXiv:1104.3630

\bibitem [NogB]{NogB} K. Noguchi. The zeta function of a finite category. arXiv:1203.6133


\end{thebibliography}
\end{document}